\documentclass{ifacconf}
\usepackage[utf8]{inputenc}
\usepackage{graphicx}      
\usepackage{natbib}        
\usepackage{bm}
\usepackage{amsmath}
\usepackage{dsfont}
\usepackage{amssymb}
\usepackage{mathrsfs}
\usepackage{pdfpages}
\usepackage{booktabs}
\usepackage{array}
\usepackage{tikz}
\usepackage{pgfplots}
\usepackage{arydshln}
\usepackage{mathtools}
\usepackage{enumitem}
\usepackage{siunitx}

\allowdisplaybreaks[3]
\makeatletter
\newcommand{\vast}{\bBigg@{5.5}}
\newcommand{\Vast}{\bBigg@{6}}

\newcounter{subeq}
\renewcommand{\thesubeq}{\theequation\alph{subeq}}
\newcommand{\newsubeqblock}{\setcounter{subeq}{0}\refstepcounter{equation}}
\newcommand{\mysubeq}{\refstepcounter{subeq}\tag{\thesubeq}}

\newtheorem{assumption}{Assumption}

\makeatletter
\DeclareFontFamily{OMX}{MnSymbolE}{}
\DeclareSymbolFont{MnLargeSymbols}{OMX}{MnSymbolE}{m}{n}
\SetSymbolFont{MnLargeSymbols}{bold}{OMX}{MnSymbolE}{b}{n}
\DeclareFontShape{OMX}{MnSymbolE}{m}{n}{
	<-6>  MnSymbolE5
	<6-7>  MnSymbolE6
	<7-8>  MnSymbolE7
	<8-9>  MnSymbolE8
	<9-10> MnSymbolE9
	<10-12> MnSymbolE10
	<12->   MnSymbolE12
}{}
\DeclareFontShape{OMX}{MnSymbolE}{b}{n}{
	<-6>  MnSymbolE-Bold5
	<6-7>  MnSymbolE-Bold6
	<7-8>  MnSymbolE-Bold7
	<8-9>  MnSymbolE-Bold8
	<9-10> MnSymbolE-Bold9
	<10-12> MnSymbolE-Bold10
	<12->   MnSymbolE-Bold12
}{}

\let\llangle\@undefined
\let\rrangle\@undefined
\DeclareMathDelimiter{\llangle}{\mathopen}%
{MnLargeSymbols}{'164}{MnLargeSymbols}{'164}
\DeclareMathDelimiter{\rrangle}{\mathclose}%
{MnLargeSymbols}{'171}{MnLargeSymbols}{'171}
\makeatother

\begin{document}
\begin{frontmatter}

\title{Distributed Frequency and Voltage Control for AC Microgrids based on Primal-Dual Gradient Dynamics\thanksref{footnoteinfo}} 

\thanks[footnoteinfo]{This work was funded by the Deutsche Forschungsgemeinschaft (DFG, German Research Foundation)---project number 360464149.}

\author[First]{Lukas K\"olsch} 
\author[First]{Katharina Wieninger} 
\author[First]{Stefan Krebs}
\author[First]{S\"oren Hohmann}

\address[First]{Institute of Control Systems, Karlsruhe Institute of Technology, 
   Karlsruhe, Germany (e-mail: lukas.koelsch@kit.edu)}

\begin{abstract}                
With the gradual transformation of power generation towards renewables, distributed energy resources are becoming more and more relevant for grid stabilization. In order to involve all participants in the joint solution of this challenging task, we propose a distributed, model-based and unifying controller for frequency and voltage regulation in AC microgrids, based on steady-state optimal control. It not only unifies frequency and voltage control, but also incorporates the classic hierarchy of primary, secondary and tertiary control layers with each closed-loop equilibrium being a minimizer of a user-defined cost function. 
By considering the individual voltage limits  as additional constraints in the corresponding optimization problem, no superordinate specification of voltage setpoints is required. Since the dynamic model of the microgrid has a port-Hamiltonian structure, stability of the overall system can be assessed using shifted passivity properties. Furthermore, we demonstrate the effectiveness of the controller and its robustness against fluctuations in active and reactive power demand by means of numerical examples.
\end{abstract}

\begin{keyword}
distributed control, optimization-based control, electric power systems, microgrids, frequency regulation, voltage regulation
\end{keyword}

\end{frontmatter}

\section{Introduction}
Our current energy system is undergoing a rapid change towards an increasing penetration by renewable energy sources. 
Large, central, conventional power plants are being successively replaced by distributed energy resources (DERs). 
If interconnected in an islanded microgrid without any superordinate instance, the individual DERs must provide stability by themselves, which has direct implications for the control strategy:
On the one hand, control actions must become increasingly faster due to the lack of rotational inertia.
On the other hand, the central superordinate authority that ensures stability is increasingly disappearing which in turn requires the control objective to be achieved jointly and decentrally by all participating resources.
According to the IEEE-PES Task Force on Microgrid Stability Definitions, Analysis, and Modeling, 
\emph{
 ``a  microgrid is  stable  if,  after  being  subjected  to  a  disturbance,  all  state variables recover to (possibly new) steady-state values which satisfy  operational  constraints  (...),  and  without  the occurrence of involuntary load shedding.''
}
(\cite{Farrokhabadi.2019}).
In terms of frequency and voltage, this entails achieving a steady state where frequencies at each node are equal and correspond to the nominal frequency and where all of the individual voltage magnitudes remain within certain limits.

Such a complex requirement can be met conveniently using optimization-based control by characterizing the desired equilibrium as an optimizer of a constrained optimization problem. 
In the existing literature, promising approaches have already been developed to enable frequency (\cite{Stegink.2017}) and voltage regulation (\cite{Magnusson2017}) in AC microgrids by this method, see \cite{Mohagheghi.2018} and \cite{Doerfler2019} for an extensive survey on current research directions on optimization-based control of future power grids.
However, these existing control approaches always pursue only one of the two objectives or implicitly assume that specific voltage setpoints are already provided by a higher level authority.

To meet the above requirements in the original and genuine sense, we present a unifying frequency and voltage controller using primal-dual gradient dynamics in this paper.
It allows real-time optimization of the connected DERs as well as conventional generators, with the primary goal of maintaining the nominal frequency and keeping all voltage magnitudes within pre-defined limits while minimizing a user-defined cost function.
Our controller is unifying in the sense that the former hierarchical division into primary, secondary and tertiary frequency and voltage control tasks is combined within a single controller and that frequency and voltage stabilization as specified above is maintained simultaneously. In particular, no communication of superordinate set points by a higher-level authority is necessary. 
The controller design is based on a nonlinear port-Hamiltonian model for AC microgrids (\cite{Koelsch2019b}) and allows for a distributed implementation.



 
The remainder of the paper is structured as follows. 
After some notational preliminaries, in section \ref{sec:sysmodel}, we briefly recall the underlying microgrid model from \cite{Koelsch2019b} and then formalize the requirements for the desired closed-loop equilibrium. 
In section \ref{sec:controller}, we derive a price-based controller for frequency and voltage regulation that meets the previously formulated requirements. Moreover, we assess shifted passivity of the subsystems and stability of the overall system.
In section \ref{sec:simulation}, we demonstrate the performance of the closed-loop system against input disturbances by means of a 12-node test network and in section \ref{sec:conclusion}, we give a brief summary of the main results.
\subsubsection{Notation}
Vector $\bm a = \mathrm{col}_i\{a_i\}=\mathrm{col}\{a_1,a_2,\ldots\}$ is a column vector of elements $a_i$, $i=1,2,\ldots$ and matrix $\bm A=\mathrm{diag}_i\{a_i\}=\mathrm{diag}\{a_1,a_2,\ldots\}$ is a (block-)diagonal matrix of elements $a_i$, $i =1,2,\ldots$. 
The $(n \times n)$-identity matrix is denoted by $\bm I_n$ and the all-ones vector is denoted by $\mathds 1$.
Positive semi-definite matrices are denoted by $\succeq 0$ and positive definite matrices or functions are denoted by $\succ 0$.
Equilibrium variables are marked with an asterisk and shifted values with respect to an equilibrium are marked with a tilde, i.e. $\widetilde{\bm x} = \bm x - \bm x^\star$. If lower and upper bounds are specified for a particular quantity $x$, these are marked with $\underline x$ and $\overline x$, respectively.
For a given $\mu \geq 0$, we define
\begin{align}
\left\llangle x\right\rrangle_\mu^+ := \left\{ \begin{array}{ll} x, &\quad \mu > 0 \vee x \geq 0 \\ 0, &\quad \text{else} \end{array}\right. . \label{plusoperator}
\end{align} 
If $\bm x$ and $\bm \mu$ are vectors of the same size $i \in \mathds N$, then \eqref{plusoperator} can be applied component-wise, i.e.
$\llangle \bm x \rrangle_{\bm \mu}^+ : \mathrm{col}_i\{\left\llangle x_i\right\rrangle_{\mu_i}^+\}$.
\section{System Model and Problem Formulation}\label{sec:sysmodel}
\subsection{Microgrid model}
\begin{figure}[t]
	\centering
	\includegraphics[width=0.85\columnwidth]{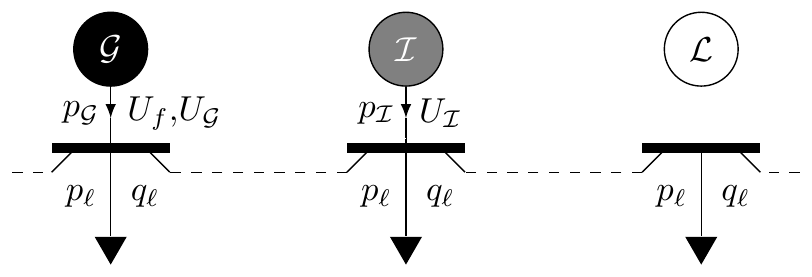}
	\caption{Schematic illustration of synchronous generator ($\mathcal G$), inverter ($\mathcal I$), and load nodes ($\mathcal L$)}	
	\label{fig:nodetypes}
\end{figure}

The microgrid is modeled by a directed graph $\mathscr G_p=(\mathcal V, \mathcal E_p)$ with $\mathcal V = \mathcal V_\mathcal G \cup \mathcal V_I \cup \mathcal V_{\mathcal L}$ being the set of $n_\mathcal G=|\mathcal V_\mathcal G|$ synchronous generator nodes , $n_\mathcal I =|\mathcal V_I|$ inverter nodes, and $n_{\mathcal L}=|\mathcal V_{\mathcal L}|$ load nodes.
These three different node types are introduced to model the connection to different microgrid participants, see Fig. \ref{fig:nodetypes}:
\begin{enumerate}[nosep]
	\item Synchronous generator nodes are connected to synchronous generators of e.g. gas or hydro turbines.
	\item Inverter nodes are connected to power electronics interfaced DERs such as photovoltaic power stations.
	\item Load nodes are connected to consumers with uncontrollable power demand or alternatively to uncontrollable power sources (modeled as negative demands).
\end{enumerate}
All nodes may be equipped with a given and uncontrollable active and reactive power demand $\bm p_\ell$ and $\bm q_\ell$ which is modelled as disturbance input. The controlled variables are active power injection $p_\mathcal G$ and excitation voltage $\bm U_f$ at generator nodes and active power injection $\bm p_\mathcal I$ and AC terminal voltage $\bm U_\mathcal I$ at inverter nodes. 
A list of all microgrid parameters used in the following can be found in Table \ref{table:PlantSymbols}.
\begin{table}
	\begin{center}
	\caption{List of Microgrid Parameters}
			\begin{tabular}{p{1.4cm}p{6.3cm}} \toprule
			\textbf{Symbol} & \textbf{Variable} \\ \midrule
			$A_i$ & positive damping coefficient \\
			$B_{ii}$ & negative of self-susceptance \\
			$B_{ij}$ & negative of susceptance of line $(i,j)$ \\
			$G_{ij}$ & negative of conductance of line $(i,j)$\\
			$L_i$ & deviation of angular momentum from nominal value $M_i\omega^{n}$ \\
			$M_i$ & moment of inertia \\ 
			$p_i$ & sending-end active power flow \\
			$p_{g,i}$ & active power generation \\
			$p_{\ell,i}$ & active power demand \\
			$q_i$ & sending-end reactive power flow \\
			$q_{\ell,i}$ & reactive power demand \\
			$U_i$ & magnitude of transient internal voltage \\
			$U_{f,i}$ & magnitude of excitation voltage \\
			$X_{d,i}-X'_{d,i}$ & d-axis synchronous minus transient reactance \\
			$\theta_i$ & bus voltage phase angle \\
			$\vartheta_{ij}$ & bus voltage angle difference\\
			$\Phi$ & overall transmission losses \\
			$\tau_{U,i}$ & open-circuit transient time constant of synchronous machine \\
			$\omega_i$ & deviation of bus frequency from nominal value\\
			\bottomrule
		\end{tabular}
		\label{table:PlantSymbols}
	\end{center}
\end{table}
The inverter interface is assumed to be equipped with an internal matching controller which allows to determine a ``virtual'' moment of inertia and a ``virtual'' damping constant for the respective nodes as shown e.g. in \cite{Jouini.2016} and \cite{Monshizadeh.}. Furthermore we make the following operating assumptions for the microgrid:
\begin{assumption} 
\begin{enumerate}
[nosep,label=\alph*)]
		\item The grid is a balanced three-phased system and the lines are represented by its one-phase $\pi$-equivalent circuits. \label{A2}
	\item The grid is operating around the nominal frequency. \label{A1}
	\item Subtransient dynamics of the synchronous generators is neglected. \label{A3}
	\item The internal matching controller of the inverters has fast dynamics compared to the gradient-based frequency and voltage controller. \label{A4}
\end{enumerate}
\end{assumption}
The physical interconnection of the nodes is represented by an incidence matrix $\bm D_p \in \mathds R^{n\times m_p}$ with $n=n_\mathcal G+n_\mathcal I+n_{\mathcal L}$ and $m_p = |\mathcal E_p|$. The incidence matrix $\bm D_p$ can be subdivided as
$\bm D_p = \mathrm{col}\{ \bm D_{p\mathcal G}, \bm D_{p\mathcal I} , \bm D_{p {\mathcal L}} \},
$
where submatrices $\bm D_{p \mathcal G}$, $\bm D_{p\mathcal I}$, and $\bm D_{p {\mathcal L}}$ correspond to the generator, inverter, and load nodes, respectively. The sending-end active and reactive power flows from node $i \in \mathcal V$ can be calculated by
the AC power flow equations [\cite{Machowski.2012}]
\begin{align}
p_i &= \sum_{j \in \mathcal N_i} B_{ij}U_iU_j \sin(\vartheta_{ij}) + G_{ii}U_i^2 && \nonumber\\
&+ \sum_{j \in \mathcal N_i}G_{ij}U_iU_j\cos(\vartheta_{ij}),&& i \in \mathcal V, \label{powerflow1}\\
q_i &= -\sum_{j \in \mathcal N_i} B_{ij}U_iU_j \cos(\vartheta_{ij}) + B_{ii}U_i^2 && \nonumber \\
&+ \sum_{j \in \mathcal N_i}G_{ij}U_iU_j\sin(\vartheta_{ij}),&& i \in \mathcal V \label{powerflow2}
\end{align}
with $\bm G$ and $\bm B$ being the conductance and susceptance matrix, respectively, and $\vartheta_{ij}=\theta_i - \theta_j$ being the voltage angle deviation between two adjacent nodes. $\mathcal N_i$ denotes the set of neighbors of $i$, i.e. $j \in \mathcal N_i$ if $(i,j) \in \mathcal E_p$ or $(j,i) \in \mathcal E_p$.

The dynamics of the individual generator, inverter and load nodes are modeled by the following descriptor system, cf. \cite{Koelsch2019b}:
\begin{align}
\dot \theta_i &= \omega_i, && i \in \mathcal V, \label{planteq1}\\
\dot L_i &= -A_i \omega_i + p_{g,i}-p_{\ell,i}-p_i, && i \in \mathcal V_\mathcal G \cup \mathcal V_\mathcal I, \label{planteq2} \\
\tau_{d,i}\dot U_i &= U_{f,i}-U_i - \frac{X_{d,i}-X'_{d,i}}{U_i} \cdot q_i, && i \in \mathcal V_\mathcal G,  \label{planteq3}\\
0 &= -A_i \omega_i -p_{\ell,i}-p_i, && i \in \mathcal V_{\mathcal L}, \label{planteq4}\\
0 &= -q_{\ell,i}-q_i, && i \in \mathcal V_{\mathcal L}. \label{planteq5}
\end{align}
Equation \eqref{planteq1} defines the relationship between nodal frequencies and angle deviations, \eqref{planteq2} describes the active power exchange between the (possibly virtual) mechanical rotor and the neighboring nodes, \eqref{planteq3} describes the transient voltage dynamics of sychronous generator nodes and \eqref{planteq4}--\eqref{planteq5} describe the active and reactive power conservation at load nodes. To obtain a state space representation, we define the plant state vector $\bm x_p$ as
\begin{align}
\bm x_p = \mathrm{col}\{\bm\vartheta,\bm L_\mathcal G,\bm L_\mathcal I,\bm U_\mathcal G,\bm\omega_{\mathcal L},\bm U_{\mathcal L}\}, \label{plant-states}
\end{align}
where for all $i \in \mathcal V$, $j \in \mathcal V_\mathcal G$, $k \in \mathcal V_\mathcal I$, $l \in \mathcal V_{\mathcal L}$:
\begin{align}
  \bm \vartheta &= \mathrm{col}_i\{\vartheta_i\}, & \bm L_\mathcal G &= \mathrm{col}_j\{L_j\}, \\
  \bm L_\mathcal I &= \mathrm{col}_k\{L_k\},   &  \bm U_\mathcal G &= \mathrm{col}_j\{U_j\},  \\
 \bm \omega_{\mathcal L} &= \mathrm{col}_l\{\omega_l\},  &  \bm U_{\mathcal L} &= \mathrm{col}_l\{U_l\}.
\end{align}
With the Hamiltonian
\begin{align}
H_p(\bm x_p) &= \frac 12 \sum_{i \in \mathcal V_\mathcal G}\left( M_i^{-1}L_i^2 + \frac{U_i^2}{X_{d,i}-X'_{d,i}}\right) \nonumber \\
& + \frac 12 \sum_{i \in \mathcal V_\mathcal I} M_i^{-1}L_i^2 \nonumber \\
&- \frac 12 \sum_{i \in \mathcal V_\mathcal G} B_{ii}U_i^2 - \sum_{(i,j) \in \mathcal E}B_{ij}U_iU_j \cos(\vartheta_{ij}) \nonumber \\
&+\frac 12 \sum_{i \in \mathcal V_{\mathcal L}} \omega_{{\mathcal L},i}^2, \label{plant-hamiltonian}
\end{align}
and the co-state $\bm z_p = \nabla H(\bm x_p)$,  
this allows to set up a port-Hamiltonian representation of \eqref{powerflow1}--\eqref{planteq5} as follows
\begin{align}
\delimitershortfall=0pt
\setlength{\dashlinegap}{2pt}
\begin{bmatrix}
\dot{\bm \vartheta} \\ \dot{\bm L}_\mathcal G \\ \dot{\bm L}_\mathcal I \\ \dot{\bm U}_\mathcal G \\ \bm 0 \\ \bm 0 \end{bmatrix}
&=\Vast[\underbrace{\begin{bmatrix}
	\bm 0 & \bm D_{p\mathcal G}^\top & \bm D_{p\mathcal I}^\top& \bm 0 &  \bm D_{p{\mathcal L}}^\top & \bm 0 \\
	-\bm D_{p\mathcal G} & \bm 0 & \bm 0 & \bm 0 & \bm 0 & \bm 0  \\
	-\bm D_{p\mathcal I} & \bm 0 & \bm 0 & \bm 0 & \bm 0 & \bm 0  \\
	\bm 0 & \bm 0 & \bm 0 & \bm 0 & \bm 0 & \bm 0  \\
	-\bm D_{p{\mathcal L}} & \bm 0 & \bm 0 & \bm 0 & \bm 0 & \bm 0 \\
	\bm 0 & \bm 0 & \bm 0 & \bm 0 & \bm 0 & \bm 0
	\end{bmatrix}}_{\bm J_p} \nonumber \\
&-\underbrace{\begin{bmatrix}
	\bm 0 & \bm 0 & \bm 0 & \bm 0 & \bm 0 & \bm 0 \\
	\bm 0 & \bm{A}_\mathcal G & \bm 0 & \bm 0 & \bm 0 & \bm 0 \\
	\bm 0 & \bm 0 & \bm{A}_\mathcal I & \bm 0 & \bm 0 & \bm 0 \\
	\bm 0 & \bm 0 &\bm 0& \bm{R}_\mathcal G & \bm 0 & \bm 0 \\
	\bm 0 & \bm 0 & \bm 0 &\bm 0& \bm{A}_{\mathcal L} & \bm 0 \\
	\bm 0 & \bm 0 & \bm 0 & \bm 0 &\bm 0& \widehat{\bm U}_{\mathcal L}
	\end{bmatrix}}_{\bm R_p}\Vast]
\bm z_p - \cdots \nonumber \\
& \cdots -
\underbrace{\begin{bmatrix}
	\bm 0 \\ \bm{\varphi}_\mathcal G \\ \bm{\varphi}_\mathcal I \\ \bm{\varrho}_\mathcal G \\ \bm{\varphi}_{\mathcal L} \\ \bm{\varrho}_{\mathcal L}
	\end{bmatrix}}_{\bm r_p}
+ 
\left[
\begin{array}{ccc:cc}
\bm 0 & \bm 0 & \bm 0 & \bm 0 & \bm 0 \\
\bm {I} & \bm 0 &\bm 0 & \bm 0 & -\bm{\widehat I}_\mathcal G \\
\bm 0 & \bm I &\bm 0 & \bm 0 & -\bm{\widehat I}_\mathcal I \\
\bm 0 &\bm 0 & \bm{\hat \tau}_U & \bm 0 & \bm 0 \\
\bm 0 & \bm 0 &\bm 0 & \bm 0 & -\bm{\widehat I}_{\mathcal L} \\
\bm 0 &\bm 0 & \bm 0 & -\bm I & \bm 0
\end{array}
\right]
\left[
\begin{array}{c}
\bm{p}_\mathcal G \\  \bm p_\mathcal I \\  \bm{U}_f \\ \hdashline[2pt/2pt] \bm{q}_\ell \\ \bm{p}_\ell  
\end{array}
\right], \label{plantPHSlossy}
\end{align}
where
\begin{alignat}{5}
\bm A_\mathcal G &= \mathrm{ diag}_i\{A_i\},&&i \in \mathcal V_\mathcal G, \\
\bm A_\mathcal I &= \mathrm{ diag}_i\{A_i\},&&i \in \mathcal V_\mathcal I, \\
\bm A_{\mathcal L} &= \mathrm{ diag}_i\{A_i\}, && i \in \mathcal V_{\mathcal L}, \\
\bm R_\mathcal G &= \mathrm{ diag}_i\left\{\left(X_{di}-X_{di}'\right)/ \tau_{U,i}\right\}, && i \in \mathcal V_\mathcal G, \\
\widehat{\bm U}_{\mathcal L} &= \mathrm{ diag}_i\{U_i\}, && i \in \mathcal V_{\mathcal L},\\
\bm \varphi_\mathcal G &= \mathrm{ col}_i\big\{G_{ii}U_i^2 + \sum_{j \in \mathcal N_i}G_{ij}U_iU_j\cos(\vartheta_{ij})\big\}, \; && i \in \mathcal V_\mathcal G, \\
\bm \varphi_\mathcal I &= \mathrm{ col}_i\big\{G_{ii}U_i^2 + \sum_{j \in \mathcal N_i}G_{ij}U_iU_j\cos(\vartheta_{ij})\big\}, && i \in \mathcal V_\mathcal I, \\
\bm \varphi_{\mathcal L} &= \mathrm{ col}_i\big\{G_{ii}U_i^2 + \sum_{j \in \mathcal N_i}G_{ij}U_iU_j\cos(\vartheta_{ij})\big\}, && i \in \mathcal V_{\mathcal L}, \\
\bm \varrho_\mathcal G &= \mathrm{ col}_i\big\{R_{g,i}\sum_{j \in \mathcal N_i} G_{ij}U_iU_j\sin(\vartheta_{ij})\big\},&&i \in \mathcal V_\mathcal G, \\
\bm \varrho_{\mathcal L} &= \mathrm{ col}_i\big\{\sum_{j \in \mathcal N_i} G_{ij}U_iU_j\sin(\vartheta_{ij})\big\},&&i \in \mathcal V_{\mathcal L}, \\
\bm{\hat \tau}_U &= \mathrm{ diag}_i\{1 / \tau_{U,i}\}, && i \in \mathcal V_\mathcal G, \\
\widehat{\bm I}_\mathcal G &= \begin{bmatrix} \bm I_{n_g \times n_g} &  \bm 0_{n_g \times n_i} & \bm 0_{n_g \times n_{\mathcal L}} \end{bmatrix},\\
\widehat{\bm I}_\mathcal I &= \begin{bmatrix} \bm 0_{n_i \times n_g} &  \bm I_{n_i \times n_i} & \bm 0_{n_i \times n_{\mathcal L}} \end{bmatrix},\\
\widehat{\bm I}_{\mathcal L} &= 
\begin{bmatrix} \bm 0_{n_{\mathcal L} \times n_g} & \bm 0_{n_{\mathcal L} \times n_i} &\bm I_{n_{\mathcal L} \times n_{\mathcal L}} \end{bmatrix}. 
\end{alignat}
The input vector in \eqref{plantPHSlossy} is composed of the control input $\bm u_p = \mathrm{col}\{\bm p_\mathcal G, \bm p_\mathcal I, \bm U_f\}$ and the disturbance input $\bm d = \mathrm{col}\{\bm q_\ell, \bm p_\ell\} $.
\subsection{Formalization of Control Objective}
The specification for zero deviation from nominal frequency $\omega^n$ and limitation of voltage magnitudes can be formalized by the following optimization problem:
\begin{align} \label{OP} \tag{OP}
\min_{\bm p_{\mathcal G}, \bm p_{\mathcal I},\bm U_f} \quad & C(\bm p_{\mathcal G}, \bm p_{\mathcal I})		&& \nonumber \\
\newsubeqblock
\mysubeq \mathrm{subject\; to}\quad & \displaystyle\Phi=\sum_{i \in \mathcal V_\mathcal G}{p_{\mathcal G,i}}+  \sum_{i \in \mathcal V_\mathcal I}{p_{\mathcal I,i}} - \sum_{i \in \mathcal V}{p_{\ell,i}},			\label{OP1.1}															  \\
\mysubeq  & \underline{\bm U}_\mathcal G \leq \bm U_\mathcal G \leq \overline{\bm U}_\mathcal G,     \label{OP1.2}         \\
\mysubeq	& \underline{\bm U}_\mathcal I \leq \bm U_\mathcal I \leq \overline{\bm U}_\mathcal I.\label{OP1.3}
\end{align}
$C(\bm p_{\mathcal G}, \bm p_{\mathcal I})$ is a user-defined, strictly convex objective function representing e.g. the electricity generation costs 
and $\Phi = \mathds 1^\top \mathrm{ col}\{\bm\varphi_\mathcal{G},\bm \varphi_\mathcal{I}, \bm\varphi_{\mathcal L}\}$ denotes the overall transmission losses. The active power balance constraint $\eqref{OP1.1}$ is necessary for zero frequency deviation, see e.g. \cite{Trip.2016}.

To enable a formulation as a distributed controller in the process of the subsequent steps, an exact reformulation of \eqref{OP} is derived where the balance constraint \eqref{OP1.1} is replaced by a sparse set of neighbor-to-neighbor balance constraints. Moreover, since $\bm U_f$ is controllable, box constraints on $\bm U_\mathcal G$ are replaced by box constraints on $\bm U_f$:
\begin{prop}
An exact reformulation of \eqref{OP} is given by
\begin{align}  \label{OP2} \tag{OP\textsuperscript{$\sharp$}}
\min_{\bm p_{\mathcal G}, \bm p_{\mathcal I}, \bm U_f, \bm \nu} \quad & C(\bm p_{\mathcal G}, \bm p_{\mathcal I})  \nonumber \\
\newsubeqblock
\mysubeq \mathrm{subject\;to}\quad & \bm D_c \bm \nu = \widehat{\bm I}_\mathcal G^\top\bm p_\mathcal G  +  \widehat{\bm I}_\mathcal I^\top\bm p_\mathcal I - \bm p_{\ell} - \bm\varphi, \label{eq-balance-komm} \\
\mysubeq	& \bm \Psi\left(\underline{\bm U}_\mathcal G\right) \leq \bm U_f \leq \bm \Psi \left(\overline{\bm U}_\mathcal G \right), \label{OP2.2} \\
\mysubeq  & \underline{\bm U}_\mathcal I \leq \bm U_\mathcal I \leq \overline{\bm U}_\mathcal I,   				  
\end{align}
where $\bm\varphi = \mathrm{col}\{\bm\varphi_\mathcal G , \bm \varphi_\mathcal I ,  \bm\varphi_{\mathcal L}\}$, $\bm D_c$ is an incidence matrix of a connected communication graph $\mathscr G_c=(\mathcal V,\mathcal E_c)$ and $\bm \Psi (\bm U_\mathcal G) = \mathrm{col}_i\{\Psi_i (U_{\mathcal G,i})\}$, $i \in \mathcal V_\mathcal G$ with 
\begin{align}
	\Psi_i (U_{\mathcal G,i})  = &U_{\mathcal G,i}\left(1+G_{ii} \left( X_{d,i}-X'_{d,i}\right)\right) \nonumber \\
	&+ \sum_{j \in \mathcal N_i} U_j \left( G_{ij} \sin \left(\vartheta_{ij}\right) - B_{ij} \cos \left(\vartheta_{ij}\right)\right).
 \label{psi}
 \end{align}
\end{prop}
\begin{pf}
	Let $\bm x_p^\diamond$ be an optimizer of \eqref{OP} and $\bm x_p^\sharp$ be an optimizer of \eqref{OP2}.
	
	To prove that  $\bm x_p^\diamond$ fulfills \eqref{eq-balance-komm}, we first recall that for each equilibrium of \eqref{plantPHSlossy} it holds that $\bm \omega = \bm 0$ (cf. Proposition 1 in \citet{Koelsch2019b}). Accordingly, \eqref{planteq2} and \eqref{planteq4} can be written as 
	\begin{align}
	\bm 0 = \widehat{\bm I}^\top_\mathcal G \bm p_\mathcal G + \widehat{\bm I}^\top_\mathcal I \bm p_\mathcal I - \bm p_\ell - \bm p \label{asdf}
	\end{align}
	with $\bm p = \mathrm{col}_i\{p_i\}$, $i \in \mathcal V$. Inserting \eqref{powerflow1} in \eqref{asdf} yields
	\begin{align}
	\bm 0 = \widehat{\bm I}^\top_\mathcal G \bm p_\mathcal G + \widehat{\bm I}^\top_\mathcal I \bm p_\mathcal I - \bm p_\ell - \bm \varphi - \bm \phi \label{dreiunddreissig}
	\end{align}
	with $\bm \phi = \mathrm{col}_i\left\{\phi_i\right\}$, $\phi_i = \sum_{j \in \mathcal N_i}G_{ij}U_iU_j \cos(\vartheta_ij)$.
	Equation \eqref{dreiunddreissig} is equivalent to \eqref{eq-balance-komm} if and only if there exists a $\bm \nu \in \mathds R^{m_c}$ with $\bm D_c \bm \nu = \bm \phi$. Since $\bm D_c$ is the incidence matrix of a connected graph, $\mathrm{rank}(\bm D_c)=n-1$ and we can delete e.g. the last row to obtain the reduced system 
	\begin{align}
	\bm D_c^\mathrm{red} \bm \nu = \bm \phi^\mathrm{red}, \label{linsysred}
	\end{align}
	 again with $\mathrm{rank}(\bm D_c^\mathrm{red})=n-1$. Note that \eqref{linsysred} is a system of $n-1$ linear equations and $m_c$ variables.
	Since $\mathscr G_c$ is assumed to be connected, $m_c$ must be greater than or equal to $n-1$. This implies that \eqref{linsysred}	is underdetermined and hence there always exists a $\bm \nu$ satisfying \eqref{eq-balance-komm}.
	
	To proove that $\bm x_p^\sharp$ fulfills \eqref{OP1.1}, left-multiply \eqref{eq-balance-komm} with $\mathds 1^\top$ which yields 
	\begin{align}
	\Phi^\sharp=\sum_{i \in \mathcal V_\mathcal G}{p_{\mathcal G,i}^\sharp}+  \sum_{i \in \mathcal V_\mathcal I}{p_{\mathcal I,i}^\sharp} - \sum_{i \in \mathcal V}{p_{\ell,i}}	.
	\end{align}
	This completes the proof of the first equivalence.
	
	To proof equivalence of \eqref{OP1.2} and \eqref{OP2.2}, let $\bm x_p^\star$ be an equilibrium of \eqref{plantPHSlossy}. From the fourth row of \eqref{plantPHSlossy} it follows that for each $i \in \mathcal V_\mathcal G$
	\begin{align}
	0 = U_{f,i}^\star - U_{\mathcal G,i}^\star - \left(X_{d,i} - X'_{d,i} \right) q_i^\star /  U_{\mathcal G,i}^\star . \label{beweis-equivalent}
	\end{align}
	Inserting \eqref{powerflow2} in \eqref{beweis-equivalent} and comparing with \eqref{psi} yields $U_{f,i}^\star = \Psi_i (U_{\mathcal G,i}^\star)$.
Since $G_{ii}\geq 0$ and $X_{d,i}-X'_{d,i} > 0$, it follows that $\nabla \Psi_i(U_{\mathcal G,i})= (1+G_{ii} ( X_{d,i}-X'_{d,i}))>0$, i.e. $\Psi_i(U_{\mathcal G,i})$ is a strictly increasing affine map, thus ${\bm U}_\mathcal G^\star  \leq \overline{\bm U}_\mathcal G  \Longleftrightarrow \bm \Psi ({\bm U}_\mathcal G^\star ) \leq \bm \Psi \left(\overline{\bm U}_\mathcal G \right)$
and 
$ \underline{\bm U}_\mathcal G \leq {\bm U}_\mathcal G^\star \Longleftrightarrow  \bm \Psi (\underline{\bm U}_\mathcal G ) \leq \bm \Psi \left({\bm U}_\mathcal G^\star \right)$. 
	This completes the proof of the second equivalence.
	
	To sum up, each $\bm x_p^\diamond$ is feasible for \eqref{OP2} and each $\bm x_p^\sharp$ is feasible for \eqref{OP}, thus \eqref{OP} and \eqref{OP2} are equivalent.
	\hfill \hfill $\blacksquare$
\end{pf}
\section{Controller Design}\label{sec:controller}
Now, for optimization problem \eqref{OP2}, a primal-dual gradient controller [\cite{JokicLazarvandenBosch2009,Stegink.2017,Stegink.}] can be applied so that together with \eqref{plantPHSlossy}, a closed-loop equilibrium is achieved which is the solution of \eqref{OP2}.
To shorten the notation, denote the vector of active power generations by $\bm p_g$, i.e. $\bm p_g = \mathrm{col}\{\bm p_\mathcal G, \bm p_\mathcal I\}$.
\subsection{Primal-Dual Gradient Controller}
\begin{prop}
	Suppose that some constraint qualification [\cite{Boyd.2015}] holds for \eqref{OP2}. Then each closed-loop equilibrium of \eqref{plantPHSlossy} together with the distributed primal-dual gradient controller
\begin{align}
\newsubeqblock \mysubeq \bm \tau_g \dot{\bm p}_{g} &= - \nabla C({\bm p}_{g})+ \widehat{\bm I}_g{\bm \lambda} + \bm u_c, \label{primal-dual-1}\\
\mysubeq \bm \tau_{\lambda}\dot{\bm \lambda}&= \bm D_{c} \bm \nu - \widehat{\bm I}_g^\top\bm p_{g}+  \bm p_{\ell} + \bm \varphi, \label{primal-dual-2} \\
\mysubeq \bm \tau_\nu \dot{\bm \nu} &= -\bm D_{c}^\top \bm \lambda, \label{primal-dual-4} \\
\newsubeqblock \mysubeq  \bm \tau_{\bm{\mu}_{\mathcal G-}}\dot{\bm{\mu}}_{\mathcal G-} &= \llangle \underline{\bm{\Psi}}_\mathcal G - \bm U_f \rrangle^+_{\bm{\mu}_{\mathcal G-}}, \label{gen-controller-begin}\\
\mysubeq \bm \tau_{\bm{\mu}_{\mathcal G+}}\dot{\bm{\mu}}_{\mathcal G+} &= \llangle  \bm U_f - \overline{\bm{\Psi}}_\mathcal G\rrangle^+_{\bm{\mu}_{\mathcal G+}}, \\
\mysubeq \bm \tau_{\bm U_\mathcal G} \dot{\bm U}_f&=  \bm \mu_{\mathcal G-} - \bm \mu_{\mathcal G+},\label{gen-controller-end} \\
 \mysubeq \bm \tau_{\bm{\mu}_{\mathcal I-}}\dot{\bm{\mu}}_{\mathcal I-} &= \llangle \underline{\bm U}_\mathcal I - \bm U_\mathcal I \rrangle^+_{\bm{\mu}_{\mathcal I-}},\label{inv-controller-begin} \\
\mysubeq \bm \tau_{\bm{\mu}_{\mathcal I+}}\dot{\bm{\mu}}_{\mathcal I+} &= \llangle  \bm U_\mathcal I - \overline{\bm{U}}_\mathcal I\rrangle^+_{\bm{\mu}_{\mathcal I+}}, \\
\bm \tau_{\bm U_\mathcal I} \dot{\bm U}_\mathcal I&=  - ( \nabla \bm{\Psi}(\bm U_\mathcal I))^\top(\bm{\mu}_{\mathcal G-} - \bm{\mu}_{\mathcal G+}) \nonumber \\
\mysubeq  & \quad - (\nabla \bm\varphi (\bm U_{\mathcal I}))^\top \bm\lambda + \bm \mu_{\mathcal I-} - \bm \mu_{\mathcal I+} \label{inv-controller-end}
\end{align}
	with $\bm u_c = - \bm \omega$, $\bm \tau >0$, $\underline{\bm{\Psi}}_\mathcal G = \bm \Psi \left(\underline{\bm U}_\mathcal G \right)$,  $\overline{\bm{\Psi}}_\mathcal G = \bm \Psi \left(\overline{\bm U}_\mathcal G \right)$ is an optimizer of \eqref{OP2}.
	\end{prop}
\begin{pf}
The Lagrangian of \eqref{OP2} is
\begin{align}
\mathscr L (\bm p_g, \bm U_\mathcal G, \bm U_\mathcal I , \bm \nu, \bm \lambda, \bm{\mu}_{\mathcal G-}, \bm{\mu}_{\mathcal G+}, \bm{\mu}_{\mathcal I-}, \bm{\mu}_{\mathcal I+}) = \nonumber \\
C(\bm p_{\mathcal G}, \bm p_{\mathcal I})
+
\bm \lambda^\top (\bm D_{c} \bm \nu - \widehat{\bm I}_g^\top\bm p_{g}+  \bm p_{\ell} + \bm \varphi) \nonumber \\
+
{\bm{\mu}}_{\mathcal G-}^\top( \underline{\bm{\Psi}}_\mathcal G -  \bm U_f ) 
+
{\bm{\mu}}_{\mathcal G+}^\top( \bm U_f - \overline{\bm{\Psi}}_\mathcal G ) \nonumber \\
+
{\bm{\mu}}_{\mathcal I-}^\top( \underline{\bm U}_\mathcal I -  \bm U_\mathcal I ) 
+
{\bm{\mu}}_{\mathcal I+}^\top( \bm U_\mathcal I - \underline{\bm{U}}_\mathcal I ).
\end{align}
The Karush-Kuhn-Tucker (KKT) conditions specifying a saddle point of $\mathscr L$ can be applied to derive a necessary condition for an optimizer of \eqref{OP2}:
	\begin{align}
	\bm 0 &= \nabla C(\bm p_g^\star) - \bm \lambda^\star,  \\
	\bm 0 &= \bm D_{c} \bm \nu^\star - \widehat{\bm I}_g^\top\bm p_{g}^\star+  \bm p_{\ell} + \bm \varphi^\star,  \\
	\bm 0 &= \bm D_{c}^\top \bm \lambda^\star,  \\
	\bm 0 &=  -  \bm{\mu}^\star_{\mathcal G-} + \bm{\mu}^\star_{\mathcal G+},  \\
	\bm 0 &= {\bm{\mu}}_{\mathcal G-}^{\star\top}( \underline{\bm{\Psi}}^\star_\mathcal G -  \bm U_f^\star ), \\
	\bm 0 &= {\bm{\mu}}_{\mathcal G+}^{\star\top}( \bm U_f^\star - \overline{\bm{\Psi}}^\star_\mathcal G ), \\
	\bm 0 &\leq  {\bm{\mu}}_{\mathcal G-}^\star, 	{\bm{\mu}}_{\mathcal G+}^\star, \\
	\bm 0 &=  (\nabla \bm{\Psi}(\bm U_\mathcal I^\star))^\top(\bm{\mu}^\star_{\mathcal G-} -  \bm{\mu}^\star_{\mathcal G+}) \nonumber \\
	& \quad + (\nabla \bm\varphi (\bm U_{\mathcal I}^\star))^\top \bm\lambda^\star
	 -  \bm{\mu}^\star_{\mathcal I-} + \bm{\mu}^\star_{\mathcal I+},  \\
\bm 0 &= 	{\bm{\mu}}_{\mathcal I-}^{\star\top}( \underline{\bm U}^\star_\mathcal I -  \bm U^\star_\mathcal I ), \\ 
	\bm 0 &= {\bm{\mu}}_{\mathcal I+}^{\star\top}( \bm U_\mathcal I^\star - \underline{\bm{U}}^\star_\mathcal I ),  \\
	\bm 0 &\leq {\bm{\mu}}_{\mathcal I-}^\star, 	{\bm{\mu}}_{\mathcal I+}^\star.
	\end{align}	
	\end{pf}
Consequently, applying the gradient method (\cite{Arrow1958}) provides \eqref{primal-dual-1}--\eqref{inv-controller-end}. \hfill $\blacksquare$

From \eqref{primal-dual-1}--\eqref{inv-controller-end} it follows that the controller is distributed, since each local controller at node $i$ only depends on local measured values as well as values of neighboring nodes $j \in \mathcal N_i$. Keeping this in mind, we now analyze convergence of the closed loop-system towards an equilibrium.
\subsection{Analysis of Shifted Passivity and Stability}
Let $\bm d^\star = \mathrm{col}\{\bm q_\ell^\star, \bm p_\ell^\star\} $ be a given, constant disturbance input vector. We assume in the following that there exists an equilibrium $\bm x^\star$ satisfying the following regularity condition:
\begin{assumption}\label{as2}
	The Hessian of $H_p(\bm x_p)$ is positive definite at steady state $\bm x_p^\star$.
\end{assumption}
This assumption can be ensured by satisfying e.g. the (relatively mild) operational condition presented in Proposition 9 of \cite{Stegink.} and originally derived in Proposition 1 of \cite{Persis.2016b}.
 
 To investigate the shifted passivity of the closed-loop system with respect to $\bm x^\star$, it is convenient to analyze plant system \eqref{plantPHSlossy}, frequency \eqref{primal-dual-1}--\eqref{primal-dual-4}, and voltage controller \eqref{gen-controller-begin}--\eqref{inv-controller-end} separately. 
\begin{prop}\label{plant-passive}
	The plant system \eqref{plantPHSlossy} with output $\bm y_p = \bm G_p^\top \bm z_p$ provides a shifted passivity property with respect to $\bm x_p^\star$ if
\begin{align}
\left[\bm z_p - \bm z_p^\star\right]^\top \left[\mathcal R(\bm x_p)- \mathcal R({\bm x_p^\star})\right] \geq 0 \label{eq-passivity}
\end{align}
with $\mathcal R(\bm x_p)= \bm R_p \bm z_p + \bm r_p$ holds. 
\end{prop}
\begin{pf}
	The plant system provides a shifted passivity property with respect to $\bm x_p^\star$ if the shifted plant Hamiltonian 
	\begin{align}
	\widetilde H_p(\widetilde{\bm x}_p):=H_p(\bm x_p)-\left( \widetilde{\bm x}_p \right)^\top \nabla H_p({\bm x_p^\star}) - H_p({\bm x_p^\star}) \label{shifted-H}
	\end{align}
	with $\widetilde{\bm x}_p = \bm x_p - \bm x_p^\star$ satisfies $\widetilde H_p(\widetilde{\bm x}_p)\succ 0$ and 
	\begin{align}
	\dot{\widetilde H}_p(\widetilde{\bm x}_p)\leq \widetilde{\bm y}_p^\top \widetilde{\bm u}_p \label{second}
	\end{align}
	with $\bm y_p = \bm G_p^\top \bm z_p$.	
	The positive definiteness condition is satisfied locally due to Assumption \ref{as2}. To investigate \eqref{second}, with the same reasoning as in \cite{Koelsch2019b}, for constant disturbance input $\bm d^\star$, \eqref{plantPHSlossy} can be expressed in the form
	\begin{align}
	\bm E\dot{\widetilde{\bm x}}_p = \bm J_p \nabla \widetilde H_p(\bm x) - \left[\mathcal R(\bm x_p)- \mathcal R(\bm x_p^\star)\right] + \bm G_p \widetilde{\bm u}_p.
	\end{align}
	With 
	$\dot{\widetilde H}_p(\widetilde{\bm x}_p)= \widetilde{\bm z}_p^\top \bm E \dot{\widetilde{\bm x}}_p$ and bearing in mind that $\bm J_p$ is skew-symmetric, this results in condition \eqref{eq-passivity}.
	\hfill $\blacksquare$ \end{pf}
Note that for lossless grids $\bm r_p = \bm 0$, and hence \eqref{eq-passivity} is always fulfilled due to the fact that $\bm R_p \succeq 0$.

Next, we examine the shifted passivity of the frequency controller \eqref{primal-dual-1}--\eqref{primal-dual-4}: 
\begin{prop} \label{prop:frequencycontroller}
The frequency controller with input $\bm u_c$ and output $\bm y_{u1}=\bm p_g$ provides a shifted passivity property if
\begin{align}
(\bm p_g - \bm p_g^\star)^\top (\nabla C(\bm p_g) - \nabla C(\bm p_g^\star)) \geq  (\bm \lambda - \bm \lambda^\star)^\top (\bm \varphi - \bm \varphi^\star). \label{bed-frequency-passive}
\end{align}
\end{prop}
\begin{pf}
By defining the frequency controller state 
$ \bm \xi_1 = \mathrm{col}\{ \bm \tau_g \bm p_g , \allowbreak \bm \tau_{\lambda}\bm \lambda, \allowbreak \bm \tau_{\nu}\bm{\nu}\}$
and the frequency controller Hamiltonian
$
H_1(\bm \xi_1)=\frac 12 \bm x_c^\top \bm \tau_c^{-1}\bm x_c
$
with
$
\bm \tau_c=\mathrm{diag}\{\bm \tau_g, \bm \tau_{\lambda},  \bm \tau_{\nu}\} \succ 0
$,
\eqref{primal-dual-1}--\eqref{primal-dual-4} can be written in port-Hamiltonian form
\begin{align}
\dot{\bm \xi}_1 &= \bm J_1 \nabla H_1(\bm \xi_1) - \bm r_1 	+ \bm G_{u1}  \bm u_c  + \bm G_{d1}   \bm p_{\ell} \label{controller1-dynamics} \\
\bm y_{u1} &= \bm G_{u1}^\top \bm \zeta_1 \label{yu1}\\
\bm y_{d1} &= \bm G_{d1}^\top \bm \zeta_1
\end{align}
where $\bm \zeta_1 = \nabla H_1(\bm \xi_1) =\mathrm{col}\{\bm p_g , \bm \lambda, \bm{\nu}\}$ and 
\begin{align}
\bm J_1 &= \begin{bmatrix} \bm 0 & \widehat{\bm I}_g & \bm 0  \\
-\widehat{\bm I}_g^\top & \bm 0 & \bm D_{c} \\
\bm 0 & -\bm D_{c}^\top& \bm 0 \end{bmatrix}, &
\bm r_1 &= \begin{bmatrix} \nabla C(\bm p_g) \\ -\bm \varphi \\  \bm 0 \end{bmatrix}, \\
\bm G_{u1} &= \begin{bmatrix} \bm I \\  \bm 0 \\ \bm 0 \end{bmatrix}, &
\bm G_{d1} &= \begin{bmatrix} \bm 0 \\  \bm I \\ \bm 0 \end{bmatrix}.
\end{align}
With the shifted controller Hamiltonian 
\begin{align}
\widetilde H_1(\widetilde{\bm \xi}_1):=H_1(\bm \xi_1)-( \widetilde{\bm \xi}_1 )^\top \nabla H_1({\bm \xi_1^\star}) - H_1({\bm \xi_1^\star}), \label{shifted-H1}
\end{align}
the shifted controller co-state $\widetilde{\bm \zeta}_1$ equals 
\begin{align}
\widetilde{\bm \zeta}_1 = \nabla \widetilde H_1(\widetilde{\bm \xi}_1) = \nabla H_1({\bm \xi}_1) - \nabla H_1({\bm \xi}_1^\star) = \bm \zeta_1 - \bm \zeta_1^\star.
\end{align}
Since the disturbance input $\bm p_\ell$ is assumed to be constant, \eqref{controller1-dynamics} can be expressed in shifted coordinates as follows
\begin{align}
\dot{\widetilde{\bm \xi}}_1 &= \bm J_1 \nabla \widetilde H_1(\widetilde{\bm \xi}_1) - [\bm r_1(\bm \xi_1) - \bm r_1(\bm \xi_1^\star)] 	+ \bm G_{u1}  \widetilde{\bm u}_c, \\
\widetilde{\bm y}_{u1} &= \bm G_{u1}^\top \widetilde{\bm \zeta_1}.
\end{align}
Due to strict convexity of $H_1({\bm \xi}_1)$, positive definiteness $\widetilde H_1(\widetilde{\bm \xi}_1)\succ 0$ is always satisfied. Hence the frequency controller is shifted passive if the time derivative 
\begin{align}
\dot{\widetilde H}_p(\widetilde{\bm x}_p) &=  ( \nabla  \widetilde H_1(\widetilde{\bm \xi}_1))^\top \dot{\widetilde{\bm \xi}}_1 \nonumber \\  
&=  \widetilde{\bm \zeta}_1^\top [\bm r_1(\bm \xi_1) - \bm r_1(\bm \xi_1^\star)] + \widetilde{\bm \zeta}_1^\top \bm G_{u1}  \widetilde{\bm u}_c 
\end{align}
fulfills $\dot{\widetilde H}_p(\widetilde{\bm x}_p)\leq \widetilde{\bm y}_{u1}^\top \widetilde{\bm u}_c $.
Bearing in mind \eqref{yu1}, this leads to $\widetilde{\bm \zeta}_1^\top [\bm r_1(\bm \xi_1) - \bm r_1(\bm \xi_1^\star)] \leq 0$, which is equivalent to condition \eqref{bed-frequency-passive}. \hfill $\blacksquare$
\end{pf}
Note that for lossless grids $\bm \varphi = \bm 0$, and hence \eqref{bed-frequency-passive} is always fulfilled due to the fact that $C(\bm p_g)$ is (strictly) convex.
 
 Since the interconnection between frequency controller and plant system is power-preserving,
 \begin{align}
 \bm u_{c1}&=-\bm y_{p1} = -\bm G_{p1}^\top \nabla H_p (\bm x_p)=-\bm \omega_g, \\
 \bm u_{p1} &= \bm y_{c1} = \bm G_{u1}^\top \nabla H_1 (\bm \xi_1) = \bm p_g,
 \end{align} 
 shifted passivity of the subsystems in terms of Propositions \ref{plant-passive} and \ref{prop:frequencycontroller} implies shifted passivity of the closed-loop system \eqref{plantPHSlossy},\eqref{primal-dual-1}--\eqref{primal-dual-4}. 
In fact, the conditions stated in Propositions \ref{plant-passive} and \ref{prop:frequencycontroller} are not necessary, since as an excess of passivity in one subsystem can compensate for the lack of passivity in the other subsystem, cf. \cite{vanderSchaft.2017}. 
Based on the previous conditions we now formulate a stability criterion for the overall system with frequency and voltage controller:
\begin{prop}
	Assume that the conditions of Propositions \ref{plant-passive} and \ref{prop:frequencycontroller} hold. For a constant input $\bm d^\star$, let $(\bm x_p^\star, \bm \xi_1^\star, \bm \xi_2^\star)$ denote an equilibrium of \eqref{plantPHSlossy},\eqref{primal-dual-1}--\eqref{inv-controller-end} and let
	\begin{align*}
	(\nabla \widetilde H_p(\widetilde{\bm U}_{\mathcal I}))^\top \dot{\widetilde{\bm U}}_{\mathcal I} 
	&<  \widetilde{\bm U}_ {\mathcal I}((\nabla \bm{\Psi}(\bm U_\mathcal I))^\top(\bm{\mu}_{\mathcal G-} - \bm{\mu}_{\mathcal G+}) \\
	& \quad - \widetilde{\bm U}_ {\mathcal I}((\nabla \bm{\Psi}(\bm U_\mathcal I^\star))^\top(\bm{\mu}_{\mathcal G-}^\star - \bm{\mu}^\star_{\mathcal G+})\\
	& \quad - \widetilde{\bm U}_ {\mathcal I}((\nabla \bm{\varphi}(\bm U_{\mathcal I}))^\top \bm \lambda - ( \nabla \bm{\varphi}(\bm U_{\mathcal I}^\star))^\top \bm \lambda^\star)
	\end{align*}
	hold. Then there exists a neighborhood $\mathcal B$ around $(\bm x_p^\star, \bm \xi_1^\star, \bm \xi_2^\star)$ such that if $(\bm x_p, \bm \xi_1, \bm \xi_2) \in \mathcal B$, then the state asymptotically converges to  $(\bm x_p^\star, \bm \xi_1^\star, \bm \xi_2^\star)$.
	\end{prop}
\begin{pf}
	With the  voltage controller state $\bm \xi_2 = \mathrm{col}\{\bm \tau_{\bm{\mu}_{\mathcal G-}}{\bm{\mu}}_{\mathcal G-},
\bm \tau_{\bm{\mu}_{\mathcal G+}}{\bm{\mu}}_{\mathcal G+} ,
\bm \tau_{\bm U_\mathcal G} {\bm U}_f, \allowbreak
\bm \tau_{\bm{\mu}_{\mathcal I-}}{\bm{\mu}}_{\mathcal I-}, \allowbreak
\bm \tau_{\bm{\mu}_{\mathcal I+}}{\bm{\mu}}_{\mathcal I+},\allowbreak 
\bm \tau_{\bm U_\mathcal I} {\bm U}_\mathcal I	 
\}$,
let $\bm \xi_2^\star$ denote an equilibrium of \eqref{gen-controller-begin}--\eqref{inv-controller-end} and define the Lyapunov  function candidate
$
\widetilde V_2(\widetilde{\bm \xi}_2) = \frac 12 \widetilde{\bm\xi}_2^\top \bm \tau_2^{-1} \widetilde{\bm\xi}_2^\top
$
where
$\bm \tau_2 = \mathrm{diag}\{\bm \tau_{\bm{\mu}_{\mathcal G-}},
\bm \tau_{\bm{\mu}_{\mathcal G+}},
\bm \tau_{\bm U_\mathcal G}, \allowbreak
\bm \tau_{\bm{\mu}_{\mathcal I-}}, \allowbreak
\bm \tau_{\bm{\mu}_{\mathcal I+}},\allowbreak 
\bm \tau_{\bm U_\mathcal I}	 
\}\succ 0$
and
 $\widetilde{\bm\xi}_2^\top = {\bm\xi}_2 - {\bm\xi}_2^\star$. This allows to set up an overall Lyapunov  function candidate $\widetilde V(\widetilde{\bm x}_p, \widetilde{\bm \xi}_1, \widetilde{\bm \xi}_2)$ for the overall closed-loop system as the sum of $\widetilde H_p$, $\widetilde H_1$, and $\widetilde V_2$:
\begin{align}
\widetilde V(\widetilde{\bm x}_p, \widetilde{\bm \xi}_1, \widetilde{\bm \xi}_2)=\widetilde H_p(\widetilde{\bm x}_p,\widetilde{\bm \xi}_2)+ \widetilde H_1(\widetilde{\bm \xi}_1) + \widetilde V_2 (\widetilde{\bm \xi}_2) \label{gesamt-lyapunov}
\end{align}
In the notation of \eqref{gesamt-lyapunov} it has been taken into account that $\bm U_\mathcal I$ is contained in the plant Hamiltonian $H_p$, see \eqref{plant-hamiltonian}. As already shown, all summands in \eqref{gesamt-lyapunov} are positive definite, thus $\widetilde V$ is also positive definite.

Respect to $\dot{\widetilde V}$, we observe that
\begin{align}
\dot{\widetilde V}(\widetilde{\bm x}_p, \widetilde{\bm \xi}_1, \widetilde{\bm \xi}_2)=\dot{\widetilde H}_p(\widetilde{\bm x}_p,\widetilde{\bm \xi}_2)+ \dot{\widetilde H}_1(\widetilde{\bm \xi}_1) + \dot{\widetilde V}_2 (\widetilde{\bm \xi}_2) \label{lyapunov-ableitung}
\end{align}
where the first summand of \eqref{lyapunov-ableitung} equals
\begin{align}
\dot{\widetilde H}_p(\widetilde{\bm x}_p,\widetilde{\bm \xi}_2)
&=
(\nabla \widetilde H_p(\widetilde{\bm x}_p))^\top \bm J_p (\nabla \widetilde H_p(\widetilde{\bm x}_p)) \nonumber \\
& \quad -  (\nabla \widetilde H_p(\widetilde{\bm x}_p))^\top \left[\mathcal R(\bm x_p)- \mathcal R({\bm x_p^\star})\right] \nonumber \\
& \quad + (\nabla \widetilde H_p(\widetilde{\bm \xi}_2))^\top \dot{\widetilde{\bm\xi}}_2 \\
&=-  (\bm z_p - \bm z_p^\star)^\top \left[\mathcal R(\bm x_p)- \mathcal R({\bm x_p^\star})\right] \nonumber \\
& \quad + (\nabla \widetilde H_p(\widetilde{\bm \xi}_2))^\top \dot{\widetilde{\bm\xi}}_2  \label{zuio}
\end{align}
due to skew-symmetry of $\bm J_p$. The second summand of \eqref{lyapunov-ableitung} is
$
\dot{\widetilde H}_1(\widetilde{\bm \xi}_1)
=
-\widetilde{\bm p}_g^\top (\nabla C(\bm p_g) - \nabla C(\bm p_g^\star))  \nonumber  +  (\bm \lambda - \bm \lambda^\star)^\top (\bm \varphi - \bm \varphi^\star).
$
The third summand of \eqref{lyapunov-ableitung} equals
\begin{align}
\dot{\widetilde V}_2(\widetilde{\bm \xi}_2) 
&=
\widetilde{\bm{\mu}}_{\mathcal G-}^\top \llangle \underline{\bm{\Psi}}_\mathcal G - \bm U_f \rrangle^+_{\bm{\mu}_{\mathcal G-}} 
+
\widetilde{\bm{\mu}}_{\mathcal G+}^\top \llangle  \bm U_f - \overline{\bm{\Psi}}_\mathcal G\rrangle^+_{\bm{\mu}_{\mathcal G+}} \nonumber \\
& \quad + 
\widetilde{\bm U}_f^\top \widetilde{\bm \mu}_{\mathcal G-} -\widetilde{\bm U}_f^\top \widetilde{\bm \mu}_{\mathcal G+} \nonumber\\
& \quad + \widetilde{\bm{\mu}}_{\mathcal I-}^\top \llangle \underline{\bm{U}}_{\mathcal I} - \bm U_{\mathcal I} \rrangle^+_{\bm{\mu}_{\mathcal I-}} 
+
\widetilde{\bm{\mu}}_{\mathcal I+}^\top \llangle  \bm U_{\mathcal I} - \overline{\bm{U}}_{\mathcal I}\rrangle^+_{\bm{\mu}_{\mathcal I+}} \nonumber\\
& \quad + 
\widetilde{\bm U}_{\mathcal I}^\top \widetilde{\bm \mu}_{\mathcal G-} -\widetilde{\bm U}_{\mathcal I}^\top \widetilde{\bm \mu}_{\mathcal G+} \label{H2-ableitung-2} \nonumber\\
& \quad + \widetilde{\bm U}_ {\mathcal I}^\top(\nabla \bm{\Psi}(\bm U_\mathcal I))^\top(\bm{\mu}_{\mathcal G-} - \bm{\mu}_{\mathcal G+}) \nonumber\\
& \quad - \widetilde{\bm U}_ {\mathcal I}^\top(\nabla \bm{\Psi}(\bm U_\mathcal I^\star))^\top(\bm{\mu}_{\mathcal G-}^\star - \bm{\mu}^\star_{\mathcal G+})\nonumber\\
& \quad - \widetilde{\bm U}_ {\mathcal I}^\top((\nabla \bm{\varphi}(\bm U_{\mathcal I}))^\top \bm \lambda - ( \nabla \bm{\varphi}(\bm U_{\mathcal I}^\star))^\top \bm \lambda^\star). 
\end{align}
In Proposition 3 of \cite{Stegink.2015}, it was shown that
$
\widetilde{\bm \mu}^\top\llangle \bm g\rrangle^+_{\bm \mu} \leq \widetilde{\bm \mu}^\top \bm g
$ and 
$
\widetilde{\bm \mu}^\top \bm g^\star  \leq \bm 0
$
holds for each convex function $\bm g = \mathrm{col}_i\{g_i\}$. Hence
\begin{align}
&\widetilde{\bm{\mu}}_{\mathcal G-}^\top \llangle \underline{\bm{\Psi}}_\mathcal G - \bm U_f \rrangle^+_{\bm{\mu}_{\mathcal G-}} 
\leq 
\widetilde{\bm{\mu}}_{\mathcal G-}^\top ( \underline{\bm{\Psi}}_\mathcal G - \bm U_f ) \nonumber \\
=
&\widetilde{\bm{\mu}}_{\mathcal G-}^\top ( \underline{\bm{\Psi}}_\mathcal G - \bm U_f^\star - \widetilde{\bm U}_f)
\leq
- \widetilde{\bm{\mu}}_{\mathcal G-}^\top \widetilde{\bm U}_f.
\end{align}
With the same procedure it can be calculated for $\mathcal G+$, $\mathcal I-$, and $\mathcal I+$ that the  first four rows of \eqref{H2-ableitung-2} are less than or equal to zero, thus
\begin{align}
\dot{\widetilde V}_2(\widetilde{\bm \xi}_2)
&\leq \widetilde{\bm U}_ {\mathcal I}^\top(\nabla \bm{\Psi}(\bm U_\mathcal I))^\top(\bm{\mu}_{\mathcal G-} - \bm{\mu}_{\mathcal G+}) \nonumber \\
& \quad - \widetilde{\bm U}_ {\mathcal I}^\top(\nabla \bm{\Psi}(\bm U_\mathcal I^\star))^\top(\bm{\mu}_{\mathcal G-}^\star - \bm{\mu}^\star_{\mathcal G+}) \nonumber \\
& \quad - \widetilde{\bm U}_ {\mathcal I}^\top(\nabla \bm{\varphi}(\bm U_{\mathcal I}))^\top \bm \lambda - (\nabla \bm{\varphi}(\bm U_{\mathcal I}^\star))^\top \bm \lambda^\star). 
\end{align}
With the assumption from Proposition \ref{plant-passive}, the first row of \eqref{zuio} is $\leq 0$ and with the assumption from Proposition \ref{prop:frequencycontroller} it also holds that $\dot{\widetilde H}_1(\widetilde{\bm \xi}_1)\leq 0$. Hence $\dot{\widetilde V}(\widetilde{\bm x}_p, \widetilde{\bm \xi}_1, \widetilde{\bm \xi}_2)\leq 0$ is fulfilled if 
\begin{align*}
(\nabla \widetilde H_p(\widetilde{\bm \xi}_2)) \dot{\widetilde{\bm\xi}}_2 
&<  \widetilde{\bm U}_ {\mathcal I}^\top(\nabla \bm{\Psi}(\bm U_\mathcal I))^\top(\bm{\mu}_{\mathcal G-} - \bm{\mu}_{\mathcal G+}) \nonumber \\
& \quad - \widetilde{\bm U}_ {\mathcal I}^\top(\nabla \bm{\Psi}(\bm U_\mathcal I^\star))^\top(\bm{\mu}_{\mathcal G-}^\star - \bm{\mu}^\star_{\mathcal G+})\nonumber\\
& \quad - \widetilde{\bm U}_ {\mathcal I}^\top(\nabla \bm{\varphi}(\bm U_{\mathcal I}))^\top \bm \lambda - (\nabla \bm{\varphi}(\bm U_{\mathcal I}^\star))^\top \bm \lambda^\star)
\end{align*}
holds. With $(\nabla \widetilde H_p(\widetilde{\bm \xi}_2))^\top \dot{\widetilde{\bm\xi}}_2
=
(\nabla \widetilde H_p(\widetilde{\bm U}_{\mathcal I}))^\top \dot{\widetilde{\bm U}}_{\mathcal I}$, the proof is complete.
\hfill $\blacksquare$	
\end{pf}
\subsection{Possible Extensions and Variations}
The notation as an optimization problem provides numerous possibilities for extending the proposed controller: 
For example, hard limitations of the active power generation $\bm p_g$ can be incorporated as additional box constraints. 
Moreover, a simplified controller with a more practical stability criterion can be formulated by assuming lossless transmission lines, i.e. $\bm \varphi = \bm 0$, and by exploiting the fact that lower and upper bounds for $\bm U_f$ are relatively insensitive with respect to $\bm U_{\mathcal I}$. Simplified' projected bounds $\underline{\bm \Psi}$ and $\overline{\bm \Psi}$ can thus be formulated by replacing all inverter voltage magnitudes $\bm U_{\mathcal I}$ in \eqref{psi} by an estimator $\widehat{\bm U}_{\mathcal I}$, e.g. $(\overline{\bm U}_{\mathcal I} + \underline{\bm U}_{\mathcal I})/2$, to obtain estimators $\underline{\widehat{\bm \Psi}}$ and $\widehat{\overline{\bm \Psi}}$ of the lower and upper bounds which are independent of $\bm U_{\mathcal I}$, leading to 
$\nabla\widehat{\bm \Psi}(\bm U_{\mathcal I}) = \bm 0 $. This gives the following simplified stability criterion:
\begin{cor}
Assume that the conditions of Proposition \ref{plant-passive} and \ref{prop:frequencycontroller} hold. For a constant input $\bm d^\star$, let $(\bm x_p^\star, \bm \xi_1^\star, \bm \xi_2^\star)$ denote an equilibrium of \eqref{plantPHSlossy},\eqref{primal-dual-1}--\eqref{inv-controller-end} with zero transmission line losses, i.e. $\bm \varphi = \bm 0$, and $\bm \Psi = \widehat{\bm \Psi}$ in \eqref{psi}. If
\begin{align*}
\sum_{i \in \mathcal V_{\mathcal I}}&\left( \mu_{\mathcal I-,i} - \mu_{\mathcal I+,i}\right) \nonumber \\
&\cdot \Bigg( \sum_{j \in \mathcal N_i} B_{ij} U_j \cos(\vartheta_{ij}) - \sum_{j \in \mathcal N_i} B_{ij} U_j^\star \cos(\vartheta_{ij}^\star)\Bigg) < 0
\end{align*}
holds, then there exists a neighborhood $\mathcal B$ around $(\bm x_p^\star, \bm \xi_1^\star, \bm \xi_2^\star)$ such that if $(\bm x_p, \bm \xi_1, \bm \xi_2) \in \mathcal B$, then the state asymptotically converges to  $(\bm x_p^\star, \bm \xi_1^\star, \bm \xi_2^\star)$.
\end{cor}
Another possible extension is to introduce an additional node type to model controllable active and reactive power input or consumption and to include reactive power limitations in the optimization problem.
\section{Simulation}\label{sec:simulation}
\begin{figure}[t]
	\centering
	\includegraphics[width=0.75\columnwidth]{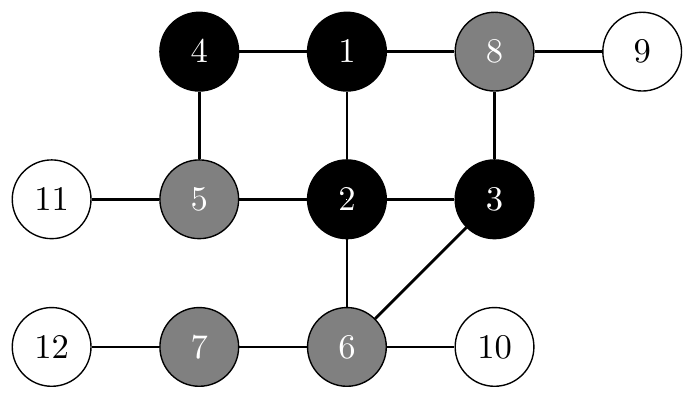}
	\caption{Network topology of exemplary microgrid.}	
	\label{fig:grid}
\end{figure}

We now verify the presented control approach by simulating an exemplary microgrid as presented in Fig. \ref{fig:grid} with $n_\mathcal G = n_\mathcal I = n_{\mathcal L} =4$.
\subsection{Parameterization}
\begin{figure}[t]
	\centering
	\includegraphics[width=\columnwidth]{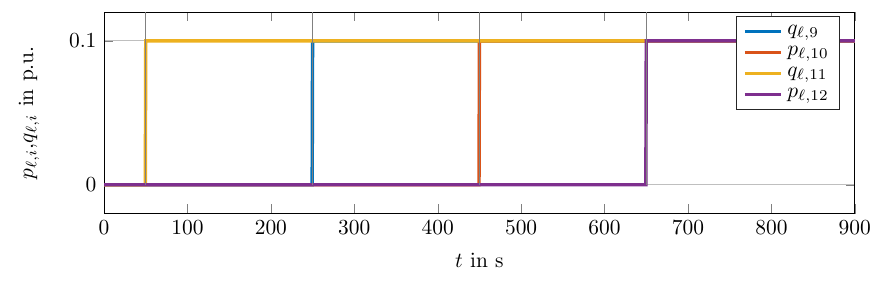}
	\caption{Stepwise increase of active and reactive power demands at load nodes.}	
	\label{fig:pl}
\end{figure}

The parameter values for all nodes, listed in Tables \ref{tab:generator}--\ref{tab:load}, are based on \cite{Trip.2016} and \cite{Koelsch2019}, whereby the ``virtual'' moments of inertia of inverter nodes were selected to be considerably smaller than the corresponding moments of inertia of generator nodes.
All values are given in p.u. ($U_{\mathrm{base}}=\SI{20}{kV},S_{\mathrm{base}}=\SI{100}{MVAr}$) except $\tau_{U,i}$, which is given in seconds. 
\begin{table}
	\begin{center}
	\caption{Parameters of Generator Nodes}
		\begin{tabular}{lllllllll}
			\toprule
			$\bm i$ 	& \textbf 1 & \textbf 2 &\textbf  3 &\textbf  4  \\
			\midrule
			$A_i$ & 1.6 & 1.22 & 1.38 & 1.42 \\
			$B_{ii}$ & -6.0567 & -8.014 & -6.6755 & -4.144\\
			$G_{ii}$ & 5.7834 & 7.86379 & 5.9209 & 4.09632\\
			$M_i$ & 26.1 & 19.9 & 22.45 & 21.1 \\ 
			$X_{d,i}$ & 0.15 & 0.19 & 0.165 & 0.1875 \\
			$X'_{d,i}$ & 0.055 & 0.045 & 0.055 & 0.056 \\
			$\tau_{U,i}$ & 6.45 & 7.68 & 7.5 & 6.5 \\
			\bottomrule
		\end{tabular}
		\label{tab:generator}
	\end{center}
\end{table}
\begin{table}
		\begin{center}
\caption{Parameters of Inverter Nodes}
		\begin{tabular}{llllllll}
			\toprule
			$\bm i$ & \textbf5 & \textbf6 & \textbf7 & \textbf8 \\
			\midrule
			$A_i$ & 1.4 & 1.3 & 1.35 & 1.45 \\
			$B_{ii}$ & -5.6611 & -8.1791 & -3.6067 & -6.1635 \\
			$G_{ii}$ & 5.77198 & 7.25506 & 3.82174& 5.74546 \\
			$M_i$ & 4.4 & 4.5 & 5.15 & 4 \\ 
			\bottomrule
		\end{tabular}
		\label{tab:inverter}
	\end{center}
\end{table}
\begin{table}
		\begin{center}
\caption{Parameters of Load Nodes}
		\begin{tabular}{lllll}
			\toprule
			$\bm i$  & \textbf 9 & \textbf{10} & \textbf{11} & \textbf{12} \\
			\midrule
			$A_i$  & 1.45 & 1.35 & 1.5 & 1.7 \\
			$B_{ii}$ & -1.716 & -2.41244 & -1.692 & -1.848 \\
			$G_{ii}$ & 2.05346 & 1.99349 & 1.83437 & 2.02776 \\
			\bottomrule
		\end{tabular}
		\label{tab:load}
	\end{center}
\end{table}
The transmission line parameters are generated choosing both line resistances $R_{ij}$ and reactances $X_{ij}$ as evenly distributed random variables around $\SI{1}{\Omega}\pm 10 \%$ and can be found in Table \ref{tab:line}. In the same fashion, the shunt capacitors at each node are set to $\SI{10}{nF}\pm 10 \%$.
The controller parameters $\bm\tau_{\bm{\mu}_{\mathcal G-}}$,$\bm\tau_{\bm{\mu}_{\mathcal G+}}$ and $\bm\tau_{\bm{U}_{\mathcal G}}$ are set to $0.01$, $\bm\tau_{\bm{U}_{\mathcal I}}$ to $10$ and all other controller parameters to $0.1$. 
The cost function is chosen to 
\begin{align}
C(\bm p_g)=\frac 12 \sum_{i \in \mathcal V_\mathcal G \cup V_\mathcal I} \frac{1}{w_i} \cdot p_{g,i}^2
\end{align}
with weighting factors $\omega_1 = 1, w_2=\num{1.1}, w_{3}=\num{1.2}$ et cetera. Since $\nabla C(p_i^\star) = \nabla C(p_j^\star)$ at steady state, this specific choice of $C(\bm p_g)$ as a weighted sum of squares leads to \emph{active power sharing}, i.e. a proportional share $p^\star_{g,i} \slash w_i  = p^\star_{g,j} \slash w_j= \text{const.}$ for all $i,j\in \mathcal V_\mathcal G \cup \mathcal V_\mathcal I$.
The voltage limits are set to $[0.98 \; 1.02]$ and the upper bound for active power generation is set to $0.6$.
We choose $\bm D_c$ to be identical to the plant incidence matrix $\bm D_p$ after it has been pointed out in \cite{Koelsch2019} that the specific choice of $\bm D_c$ has little influence on the convergence speed to the desired equilibrium.
The initial values of disturbance input vector $\bm d = \mathrm{col}\{\bm p_{\ell}, \bm q_{\ell}\}$ and state vector $\bm x = \mathrm{col}\{\bm x_p, \bm \xi_1, \bm \xi_2\}$ are chosen such that the closed-loop system starts in synchronous mode with $\bm \omega ( t=0) = \bm 0$ and such that a number of voltage magnitudes are already close to or at their limits of 0.98 or 1.02.
\begin{table}[t]
	\centering
	\caption{Parameters of Transmission Lines}
	\begin{tabular}{lll}
		\toprule
		$\bm{(i,j)}$ & $\bm{B_{ij}}$ &  $\bm{G_{ij}}$\\
		\midrule
		$(1,2)$ & 1.905 &-1.9167\\
		$(1,4)$ & 1.976 &-2.04\\
		$(1,8)$ & 2.176 &-1.19178\\
		$(2,3)$ & 2.352 &-2.3256\\
		$(2,5)$ & 1.966 &-1.66 \\
		$(2,6)$ & 1.8012 &-1.8444\\
		$(3,6)$ & 2.052 &-1.7396\\
		\bottomrule
	\end{tabular}
	\hspace{0.2cm}
	\begin{tabular}{lll}
		\toprule
		$\bm{(i,j)}$ & $\bm{B_{ij}}$ &  $\bm{G_{ij}}$\\
		\midrule
		$(3,8)$ & 2.2716 &-1.7394\\
		$(4,5)$ & 2.168 &-2.016\\
		$(5,11)$ & 1.692 &-1.7984\\
		$(6,7)$ & 1.7588 &-1.7588\\
		$(6,10)$ & 2.41244&-1.9544 \\
		$(7,12)$ & 1.848 &-1.988\\
		$(8,9)$ & 1.716 &-2.0132 \\
		\bottomrule
	\end{tabular}
	\hspace{0.2cm}
	\label{tab:line}
\end{table}
At regular intervals of \SI{200}{\second}, a step of +\num{0.1}\,p.u. active or reactive power demand is assigned sequentially at load nodes, as shown in Fig. \ref{fig:pl}. 
The simulations are carried out in Wolfram Mathematica 12.0.
\begin{figure}[t!]
	\centering
	\includegraphics[width=\columnwidth]{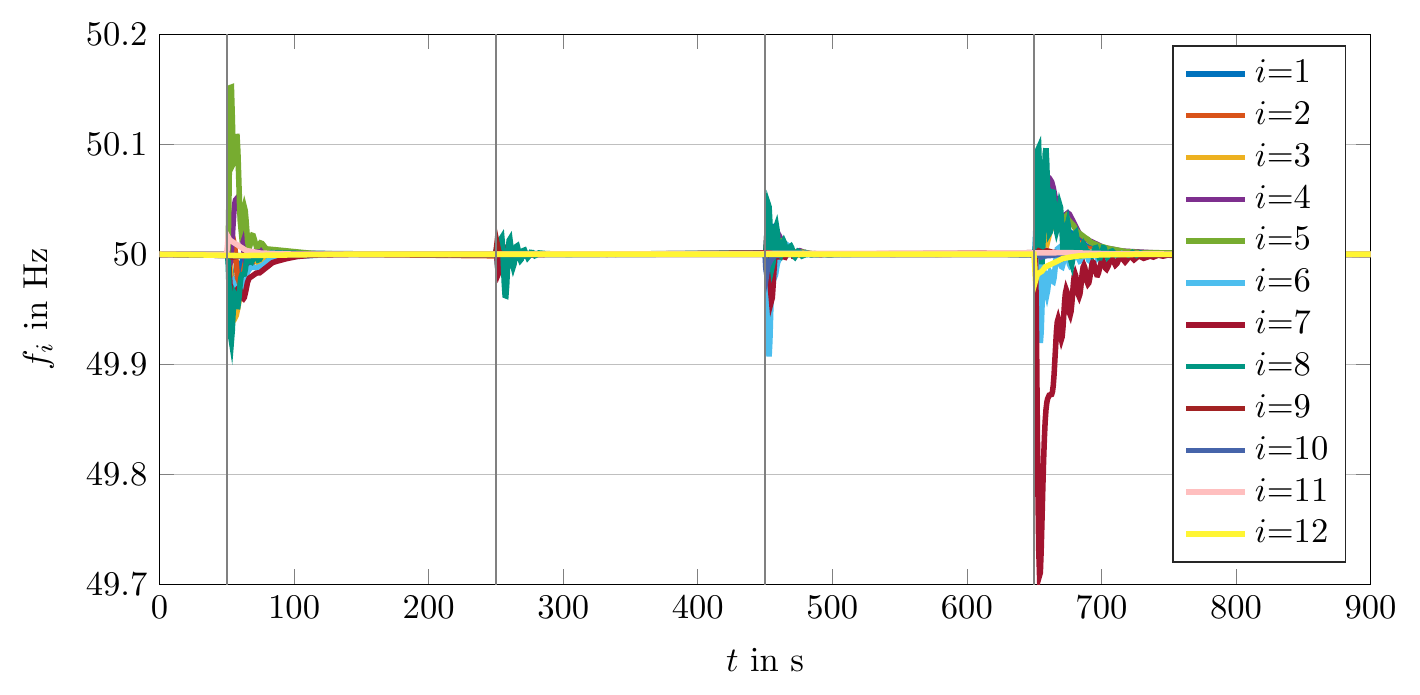}
	\includegraphics[width=\columnwidth]{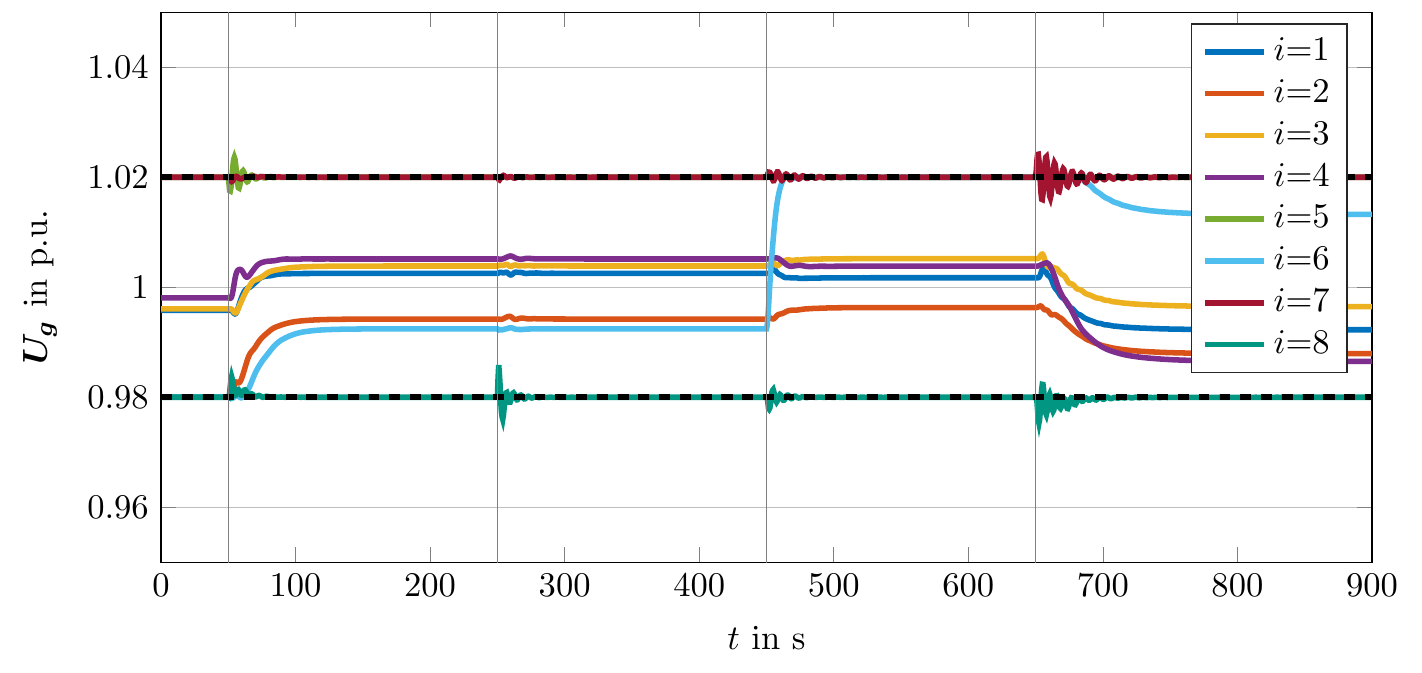}
	\includegraphics[width=\columnwidth]{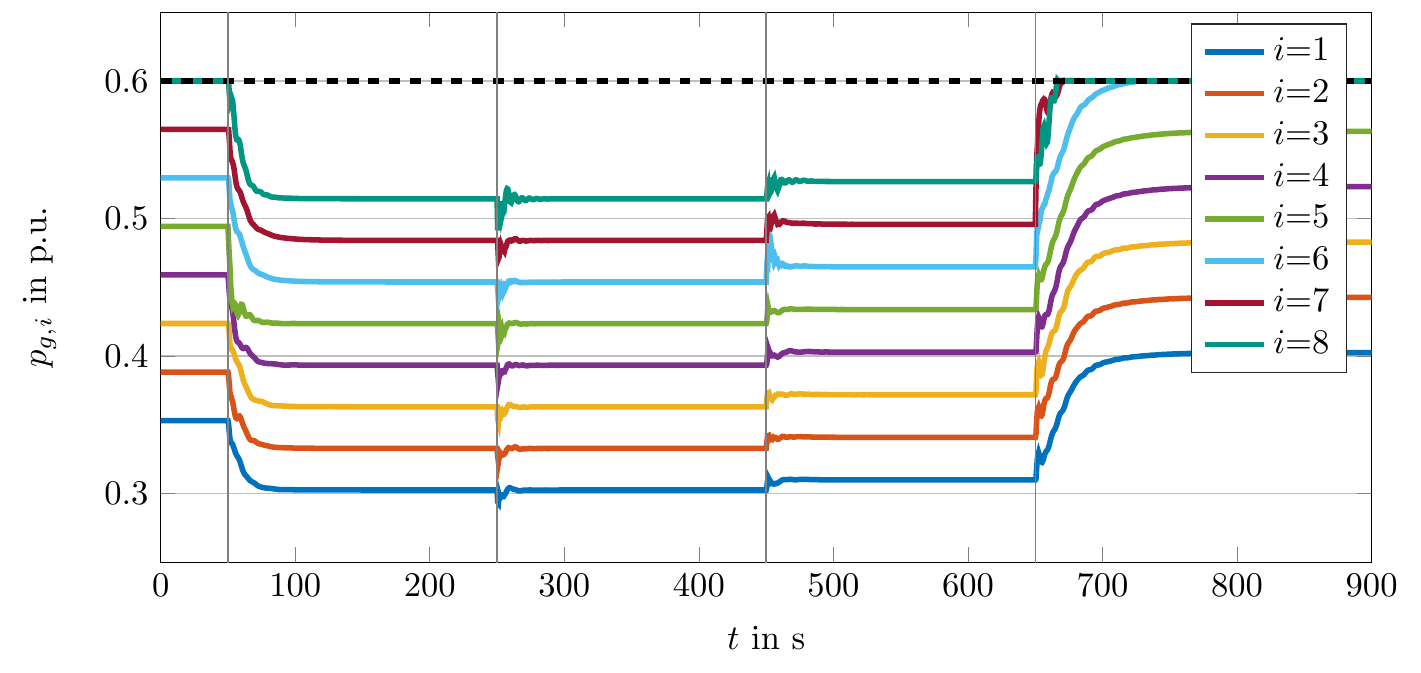}
	\caption{Evolution of frequencies (a), voltage magnitudes (b), and active power generation (c) after step increase}	
	\label{fig:Res}
\end{figure}
\subsection{Numerical Results}
Fig. \ref{fig:Res}(a) and \ref{fig:Res}(b) show the node frequencies and the voltage magnitudes, respectively. It can be seen that in steady state all frequencies are synchronized to \SI{50}{Hz} and the voltage limits of 0.98 and 1.02 are kept. 
The overshoots are around \SI{0,3}{Hz} for the node frequencies and \num{0,005}\,{p.u.} for the voltage amplitudes. 
Fig.  \ref{fig:Res}(c) shows the corresponding active power generations $\bm p_g$ at generator and inverter nodes. Remarkably, the individual power injections $p_{g,i}$ are equidistant from each other at steady state, regardless of the total generation, thus active power sharing is given.
In addition, it can be stated after the jump at \SI{650}{s} that the specified maximum limit of 0.6 {p.u.} for active power generation is not exceeded. Although this limit is reached for nodes 7 and 8, the remaining nodes perform active power sharing without being affected. 

\begin{figure}[t]
	\centering
	\includegraphics[width=\columnwidth]{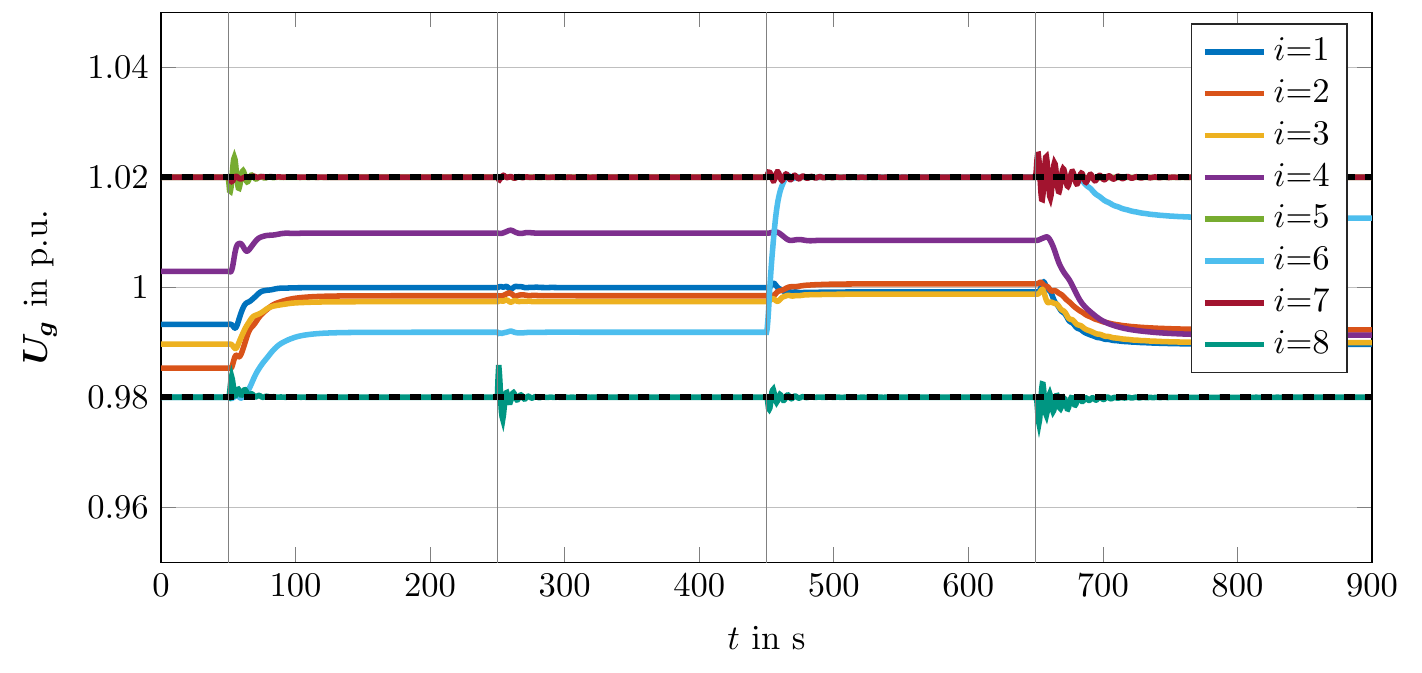}
	\caption{Evolution of voltage magnitudes if $\bm \Psi = \widehat{\bm \Psi}$.}	
	\label{fig:Uapprox}
\end{figure}
Fig. \ref{fig:Uapprox} shows the voltage magnitudes for the same scenario as above in case the simplified upper and lower limits $\widehat{\overline{\bm \Psi}}$ and $\widehat{\underline{\bm \Psi}}$ are used. 
Again, overshoots of about \num{0,005}\,{p.u.} can be detected. The progression of voltage magnitudes over time is slightly different, but comparable to Fig. \ref{fig:Res}(b).
In particular, the voltage limits are not exceeded at steady state.
\section{Conclusion}\label{sec:conclusion}
In this paper, we presented a model-based frequency and voltage controller for AC microgrids ensuring optimality with regard to a user-defined cost function.
For this purpose, the controller continuously tracks the KKT point of a constrained optimization problem using the gradient method. The underlying nonlinear microgrid model consists of a mixture of conventional synchronous generators, power electronics interfaced sources and uncontrollable loads.
The user-defined cost function can be specified in such a way that e.g. active power sharing is achieved.
Moreover, the resulting controller equations exhibit a distributed neighbor-to-neighbor communication structure.




\bibliography{IFAC_WC_Quellen}             

\appendix
\end{document}